\documentclass[a4paper,11pt]{article}

\usepackage{amsthm}
\usepackage{amsmath}
\usepackage{amssymb}
\usepackage{bbm}
\usepackage{mathtools}
\usepackage{hyperref}
\usepackage{mathrsfs}
\usepackage{xcolor}
\usepackage{MnSymbol}
\usepackage{fancyhdr}
\usepackage[nottoc, notlot, notlof]{tocbibind}
\usepackage{pdfpages}

\newcommand{\st}{\tilde{s}}

\newcommand{\dist}{\mathrm{dist}}

\newcommand{\x}{\times}
\newcommand{\p}{\partial}

\newcommand{\curl}{\mathrm{curl}}

\renewcommand{\d}{\mathrm{d}}

\newcommand{\bigO}[1]{\mathcal{O}\left(#1\right)}

\newcommand{\Lap}{\mathrm{\Delta}}

\newcommand{\eqnb}{\begin{equation}}
\newcommand{\eqnbs}{\begin{equation*}}
\newcommand{\eqnbsa}{\begin{equation*}\begin{aligned}}
\newcommand{\eqnba}{\begin{equation}\begin{aligned}}
\newcommand{\eqnbl}[1]{\begin{equation}\label{#1}}
\newcommand{\eqnbal}[1]{\begin{equation}\label{#1}\begin{aligned}}
\newcommand{\eqnes}{\end{equation*}}
\newcommand{\eqne}{\end{equation}}
\newcommand{\eqnesa}{\end{aligned}\end{equation*}}
\newcommand{\eqnea}{\end{aligned}\end{equation}}

\newcommand{\re}[1]{(\ref{#1})}
\newcommand{\comment}[1]{}

\newcommand{\RR}{\mathbb{R}}
\newcommand{\TT}{\mathbb{T}}

\newcommand{\ZZ}{\mathbb{Z}}

\newcommand{\cH}{\mathcal{H}}
\newcommand{\cC}{\mathcal{C}}

\newtheorem{theorem}{Theorem}

{\bf}{\it}

\title{Asymptotics for vortex filaments using a modified Biot-Savart kernel}
\author{Benjamin C.\ Pooley$^1$  and Jos\'e L.\ Rodrigo$^2$}
\date{\it$^{1,2}$Mathematics Institute, University of Warwick,\\ Coventry, CV4 7AL}
\begin{document}
\maketitle
\abstract{
We consider a family of approximations to the Euler equations obtained by adding $(-\Lap)^{-\alpha/2}$ to the non-locality in the Biot-Savart kernel together with a mollification (with parameter $\varepsilon$).  We consider the evolution of a thin vortex tube. We show that  the velocity on the filament (core of the tube) in the limit as $\varepsilon\to 0$ is  given $\frac{C(\alpha,t)}{\alpha} \kappa B + \mathcal O(1)$ where $\kappa$ and $B$ are the curvature and binormal of the curve, and $C$, $C^{-1}$ are uniformly bounded. }
\section{Introduction}
The Euler equations model the evolution of an inviscid incompressible Newtonian fluid. In 3-dimensions they can be formulated in terms of the vorticity $\omega:[0,T]\x\RR^3\to \RR^3$, which is the curl of the fluid velocity $\omega=\nabla\x u$:

\begin{equation}\label{eq_vort}
\omega_t+(u\cdot\nabla)\omega=(\omega\cdot\nabla)u,
\end{equation}
\begin{equation}\label{eq_BS}
u=\curl^{-1}\omega,
\end{equation}
where $\curl^{-1}$ in $\RR^3$ is given  by the Biot--Savart operator on $\RR^3$: 
\begin{equation}\label{eq_BS}
\curl^{-1}\omega(x)\coloneqq \frac{-1}{4\pi}\int_{\RR^3}\frac{x-y}{|x-y|^3}\x\d\omega(y),
\end{equation}
where $\d\omega(y)=\omega(y)\d y$ if $\omega$ is a function. This is a bounded operator $L^p\to L^q$ for example if $p\in (1,3)$ and $1/p = 1/q + 1/3$ (see \cite{JCR_NSE_Book}, for example).

Describing the evolution of $\omega$ in the case that it is initially highly concentrated around a  {\it filament},  and approximately tangential to the curve is a long-standing problem. In fact it dates back at least as far as the 1860s. See for example the paper by Helmholtz \cite{Helmholtz1867} and the letter by Kelvin appended to it. In the early 20th century, the asymptotic {\it local induction approximation} was  developed by Da Rios (see \cite{DaRios1906},\cite{Ricca1996}, and \cite{RiccaNature_1991}) and much later by Arms and Hama \cite{ArmsHama_1965}, among other authors. A more recent treatment can be found in \cite{TingBook2007}.

For sufficiently regular flows integral curves of the vorticity field are advected by the velocity. Assuming the same to be true for an isolated vortex line, consider a time dependent curve $\cC(t)$ that is advected by the velocity $u$ from \re{eq_BS} where $\omega$ is the $\cH^1$ measure on $\cC$ multiplied by the tangent. The local induction approximation states that an advected point $p(t)\in \cC(t)$ satisfies
\begin{equation}\label{eqLIA}
\frac{\d }{\d t}p(-t/\log(\varepsilon))\big|_{t=0}\asymp\kappa B
\end{equation}
 as $\varepsilon\to 0$, where $\kappa$ and $B$ are the curvature and binormal to $\cC$ at $p$. 
 
 Recently Jerrard and Seis \cite{JerrardSeis_2017} found new estimates of the difference between the evolution of a weak solution of Euler in 3-dimensions with approximately-filamentary vorticity and the binormal curvature flow. This was based on the notion of weak binormal curvature flows developed in \cite{JerrardSmets2015}. 
 
 The evolution of a curve according to the binormal curvature flow is a subject for study in its own right and is related to a nonlinear Schr\"odinger equation via the Hasimoto transform \cite{Hasimoto1972}. For recent work on this topic,  see \cite{HozVega2018}, \cite{delaHoz_Vega_2014}, \cite{BanicaVega2015}, or \cite{Vega2015}, for example.
 
Although this paper is concerned with three-dimenional flows, it is worth noting that the analogous problem for two-dimensional Euler concerns the evolution of systems of point vortices. Analysis of such systems has to date proved more fruitful than that of their three-dimensional counterparts. The classical {\it vortex model} states that the evolution of a  collection of vortex points $\{x_k(t)\}$ is such that each advected by the velocity corresponding to others:
 \[
 \dot x_k = \sum_{n\neq k}\curl^{-1}\delta_{x_n}.
 \]
This has been properly justified by Marchioro and Pulvirenti \cite{MarchioroPulvirenti1993} (see also \cite{Marchioro_Pulvirenti_1994, MarchioroPulvirenti83}), who showed that until such a time that two of the points collide under this evolution, the vorticity of a solution of the Euler equations is concentrated at the points $x_k(t)$, if initially it is sufficiently concentrated at the points $x_k(0)$. 

Recently Davila et.\ al.\ \cite{Davila_delPino_etal2018} have found a way to construct solutions of the Euler equations in two dimensions (including in bounded domains) with vorticity uniformly approximating a sum of desingularised $\delta$-distributions moving according to the Kirchoff--Routh law, which generalises the simple vortex model above for bounded domains. Whereas the approach of Marchioro and Pulvirenti was largely concerned with controlling the support of vorticity,  solutions constructed using these explicit desingularisations admit finer information about the vorticity inside the core. 

For a concise survey of some of these topics, see  \cite{BanicaMiot2013}.

In order to avoid the need for time rescaling of the form \re{eqLIA}, we study the model of the Euler equations in $\RR^3$ described below, in the case that $\omega$ is a vector-valued measure given by the tangent to a closed simple smooth curve $\mathcal{C}$, parametrised by $\gamma\in C^\infty(\TT^1;\RR^3)$, with non-vanishing tangent $\gamma'$. Here $\TT^1$ denotes the torus $\RR/\ZZ$. We want to replace the velocity $u$ given by $\curl^{-1} \omega$ in  \re{eq_BS} by:

\begin{equation}\label{eq_alph_eps_BS}
u=u_\varepsilon^\alpha=\mathcal{J}_\varepsilon\curl^{-1}\Lambda^{-\alpha}\omega,
\end{equation}
where $\mathcal{J}_\varepsilon$ denotes mollification by some fixed three-dimensional mollifier $\eta_\varepsilon=\varepsilon^{-3}\eta(x/\varepsilon)$, $\Lambda=(-\Lap)^{1/2}$, $\varepsilon>0$, and $\alpha\in(0,1/2)$.

Using the notation  $\cH^1_\cC\coloneqq\cH^1\lefthalfcup \cC$ for the 1-Hausdorff measure restricted to $\cC$, we can state the main result of this paper as follows.
\begin{theorem}\label{thm1}
If $\cC$ is a smooth curve with smooth parametrisation $\gamma$, and $\omega=\gamma'\cdot \cH^1_\cC$  then for $u_\varepsilon^\alpha$ defined by \re{eq_alph_eps_BS} satisfies
\[
\lim_{\varepsilon\to 0}u_\varepsilon(\gamma(\tau))=\frac{C(\alpha,\tau)}{\alpha}\kappa(\tau)B(\tau)+ w(\alpha,\tau),
\]
where $C>0$ is bounded above and below independent of $(\alpha,\tau)$, $|w|$ is bounded independent of $(\alpha,\tau)$, and $\kappa$ and $B$ denote the curvature and binormal to $\cC$ at $\gamma(\tau)$.
 \end{theorem}
 
 The mollification appearing in \re{eq_alph_eps_BS} effectively removes a singularity of order $\dist(x,\mathcal{C})^{\alpha -1}$ as $x$ approaches the curve $\mathcal{C}$. This corresponds to a component of the velocity of the order $\dist(x,\mathcal{C})^{-1}$, that also appears in the analysis of the Euler equations. In the classical case, that term is not usually considered to play a role in the evolution of the curve.  The remaining terms are much more remarkable. Indeed, for $\alpha>0$,  Theorem \ref{thm1} implies that after removing the first singularity, the velocity of the curve is finite and moreover, it is only the magnitude that depends on $\alpha$, to leading order. This is in stark contrast to vortices in the classical system, where the binormal term is also singular, of order $|\log(\dist(x,\mathcal{C}))|$, which warrants the time-rescaling seen in \ref{eqLIA}.
 
 We remark that Berselli and Gubinelli \cite{BerselliGubinelli2007} have proved well-posedness results for filaments $\gamma$ evolving  under velocities of the form
 \[
 u(t,x)=\int_{0}^1\nabla \phi(x-\gamma(t,s))\x\p_s\gamma(t,s)\,\d s,
 \]
 where $\phi$ is even with integrable and non-negative Fourier transform $\hat\phi$ and $\int_{\RR^3} (1+|\xi|^2)^2\hat\phi(\xi)\,\d \xi <\infty$.  Such examples include the so-called Rosenhead approximation \cite{Rosenhead1930} ($\phi(x)=c(|x|^2+\mu^2)^{-1/2}$ for some $\mu\neq0$), but not the system
 \[
 u=\curl^{-1}\Lambda^{-\alpha}\omega.
 \]
 
The rest of the paper is dedicated to the proof of Theorem \ref{thm1}.

\section{Taylor expanding kernel of $\curl^{-1}\Lambda^{-\alpha}$}

In this section we calculate $v=\curl^{-1}\Lambda^{-\alpha}\omega$ in the vicinity of $\mathcal{C}$

Let $\gamma$ be a  smooth, closed and simple curve; it admits a {\it security radius}
$r_s>0$, such that for all $x\in B_{r_s}(\mathcal{C})$ there exists a unique $\tau\in\TT^1$ for which $|x-\gamma(\tau)|<r_s$ and $(x-\gamma(\tau))\cdot\gamma'(\tau)=0$. In this case we say that $x$ is {\it within the security radius at $\tau$}.

To simplify the notation let us assume that $\gamma(0)=0$ and consider $x\neq 0$ within the security radius at $0$. Calculating $u$ within the security radius at other points can be achieved by changing variables. The modified Biot-Savart law corresponding to $\curl^{-1}\Lambda^{-\alpha}$ is
\begin{equation}\label{eq_aBS_ker}
v(x)=c_\alpha\int_{\RR^3} \frac{x-y}{|x-y|^{3-\alpha}}\x\d\omega(y),
\end{equation} 
where $c_\alpha$ is bounded independent of $\alpha$ and we will omit it from the following calculations.

If $\omega=\gamma'\cdot \cH^1_\cC$,  \re{eq_aBS_ker} becomes
\begin{equation}\label{aux00001}
v(x)=\int_{\TT^1}\frac{x-\gamma(\st)}{|x-\gamma(\st)|^{3-\alpha}}\x\gamma'(\st)|\gamma'(\st)|\ \d \st.
\end{equation}

We now calculate the Taylor expansion of the following term in the integrand
\[
\frac{x-\gamma(\st)}{|x-\gamma(\st)|^{3-\alpha}}\x\gamma'(\st),
\]
with respect to $\st$,  about $\st=0$. For definiteness, we consider $\st$ in the fundamental domain $[-1/2,1/2)$ of $\TT^1$. Let $R=R(x,\st)$ be defined by
\[
R\coloneqq \sqrt{|x|^2+|\gamma'(0)|^2\st^2},
\]
then for $x$ in the security radius (i.e.\ $x\in B_{r_0}(0)$)
\begin{multline*}
|x-\gamma(\st)|^2=R^2\left(1-\frac{\st^2x\cdot\gamma''(0)}{R^2}-\frac{\st^3(x\cdot\gamma^{(3)}(0)-3\gamma'(0)\cdot\gamma''(0))}{3R^2}\right.\\
\left.+\frac{\st^4(3|\gamma''(0)|^2+4\gamma'(0)\cdot\gamma^{(3)}(0))}{12R^2}+\bigO{\frac{|x|\st^4}{R^2},\frac{\st^5}{R^2}}\right).
\end{multline*}
Hence for the denominator we obtain
\begin{multline*}
|x-\gamma(\st)|^{\alpha-3}=R^{\alpha-3}\left(1+\frac{3-\alpha}{2}\left(\frac{\st^2x\cdot\gamma''(0)}{R^2}\right.\right.\\
\left.\left.+\frac{\st^3(x\cdot\gamma^{(3)}(0)-3\gamma'(0)\cdot\gamma''(0))}{3R^2}-\frac{\st^4(3|\gamma''(0)|^2+4\gamma'(0)\cdot\gamma^{(3)}(0))}{12R^2}\right)\right.\\
\left.+\bigO{\frac{|x|\st^4}{R^2},\frac{\st^5}{R^2},\frac{\st^4|x|^2}{R^4},\frac{\st^8}{R^4}}\right).
\end{multline*}
Now, for the numerator 
\begin{multline*}
(x-\gamma(\st))\x\gamma'(\st)= x\x\gamma'(0) + \st x\x\gamma''(0)\\
-\frac{\st^2}{2}\gamma'(0)\x\gamma''(0)+\bigO{\st^3,|x|\st^2}.
\end{multline*}
Therefore
\begin{multline}\label{eq_expn_unmol}
\frac{x-\gamma(\st)}{|x-\gamma(\st)|^{3-\alpha}}\x\gamma'(\st)\\= \frac{1}{R^{3-\alpha}}\left( x\x\gamma'(0) + \st x\x\gamma''(0)\vphantom{\frac{\st^2}{2}}-\frac{\st^2}{2}\gamma'(0)\x\gamma''(0)\right)\\
+\frac{3-\alpha}{2R^{5-\alpha}}\left(\st^2x\cdot\gamma''(0)-\st^3\gamma'(0)\cdot\gamma''(0)\right)x\x\gamma'(0)\\
+\bigO{\frac{\st^3}{R^{3-\alpha}},\frac{|x|\st^2}{R^{3-\alpha}},\frac{|x|^2\st^3}{R^{5-\alpha}},\frac{|x|\st^4}{R^{5-\alpha}},\frac{\st^5}{R^{5-\alpha}},\frac{\st^4|x|^3}{R^{7-\alpha}},\frac{\st^6|x|^2}{R^{7-\alpha}},\frac{\st^8}{R^{7-\alpha}}}.
\end{multline}
\section{Mollifying the expansion}

Since we are interested in a velocity field given by \eqref{eq_alph_eps_BS}, we need to consider the contribution to $u_\varepsilon$ of each term in \re{eq_expn_unmol}, via \eqref{aux00001}. 

We want to consider the mollified velocity at a point on the curve. To apply the mollification at a given $\tau\in\TT^1$ we fix a smooth orthonormal frame $n_1(\tau), n_2(\tau)$ spanning $\gamma'(\tau)^\perp$ for $\tau\in\TT^1$, this induces the following change of coordinates within the security radius of the curve $\gamma$:
\[
\Psi(\tau,y_1,y_2)=\gamma(\tau)+y_1n_1(\tau)+y_2n_2(\tau).
\]
We also define
\[
\overline{\Psi}(\tau,y)\coloneqq y_1n_1(\tau)+y_2n_2(\tau),
\]
for the  component orthogonal to the curve, and note that
\[
\left|\det\nabla\Psi\right|=|\gamma'(\tau)|+\bigO{|y|}.
\] 
\subsection{Leading term}
A careful argument using the mollification will allow us to reduce the order of the first term. Fix $\tau\in\TT^1$ and $\varepsilon\in(0,r_{\tau})$, $u_\varepsilon$ on the curve $\gamma$ is given by
\begin{multline*}
u_\varepsilon(\gamma(\tau))=\int_{\TT^1}\int_{\RR^2}\eta_\varepsilon(\gamma(\tau)-\Psi(s,y))\left|\det\nabla\Psi\right|\cdot\\
\int_{\TT^1}\frac{\overline{\Psi}(s,y)-\tilde\gamma(\st)}{|\overline{\Psi}(s,y)-\tilde\gamma(\st)|^{3-\alpha}}\x\tilde\gamma'(\st)|\gamma'(\st+s)|\,\d \st\,\d y\,\d s,
\end{multline*}
by periodicity, where $\tilde\gamma(\st)=\gamma(\st+s)-\gamma(s)$. Note that $\tilde \gamma$ is a curve for which the expansion from the previous section holds, as we assumed at the time that $\gamma(0)=0$.

Thus the first term in the expansion \re{eq_expn_unmol} contributes the following 
\begin{multline*}
\int_{\TT^1}\int_{\gamma'(s)^\perp}\eta_\varepsilon(\gamma(\tau)-\gamma(s)-z)[|\gamma'(s)|+\bigO{|z|}]\cdot\\
\int_{-1/2}^{1/2}\frac{z\x\gamma'(s)}{R^{3-\alpha}}|\gamma'(\st+s)|\,\d \st\,\d z\,\d s
\end{multline*}
where now $R=\sqrt{|z|^2+|\gamma'(s)|^2\st^2}$.

Now to make use of the anti-symmetry of $z\x\gamma'(s)$ in $\gamma'(s)^\perp$ we decompose the mollifier as follows:
\begin{multline*}
\eta_\varepsilon(\gamma(\tau)-\gamma(s)-z)=\eta_\varepsilon((\tau-s)\gamma'(s) - z ) \\
+\bigO{\varepsilon^{-2}}\chi_{\{(s,z)\colon |(\tau-s)\gamma'(s) - z|\leq C_1\varepsilon \}}
\end{multline*}
by the Mean Value Theorem, for some $C_1>0$. Indeed, since $\gamma$ is a smooth parametrisation of a smooth curve, there exists $C'>0$ such that for sufficiently small $\varepsilon>0$, 
\[\Psi^{-1}(B_\varepsilon(\gamma(\tau)))\subset (\tau-C'\varepsilon,\tau+C'\varepsilon)\x B_{C'\varepsilon}(0)\]
 for all $\tau\in\TT^1$. In which case, both $\eta_\varepsilon(\gamma(\tau)-\gamma(s)-z)$, and $\eta_\varepsilon((\tau-s)\gamma'(s) - z )$ vanish if $|(\tau-s)\gamma'(s) - z|> \varepsilon\max(\sqrt{2}C',1)$. To save notation, define 
 \[\Sigma(\tau,\varepsilon)\coloneqq\{(s,z)\in\TT^1\x\gamma'(s)^\perp\colon |(\tau-s)\gamma'(s) - z|\leq C_1\varepsilon \}.\]

We can also absorb the $\bigO{|z|}$ part of the determinant into the error term:
\begin{multline*}
\eta_\varepsilon(\gamma(\tau)-\gamma(s)-z)[|\gamma'(s)|+\bigO{\varepsilon}]\\ 
=\eta_\varepsilon((\tau-s)\gamma'(s) - z )|\gamma'(s)| +\bigO{\varepsilon^{-2}}\chi_{\Sigma(\tau,\varepsilon)}.
\end{multline*}

By the oddness with respect to $z$, the contribution of the following term vanishes (since $|(\tau-x)\gamma'(x) - z|=|(\tau-x)\gamma'(x) + z|$, and $\eta$ is chosen to be a radial function):
\[
\int_{\gamma'(s)^\perp}\eta_\varepsilon((\tau-s)\gamma'(s)-z)|\gamma'(s)|\int_{-1/2}^{1/2}\frac{z\x\gamma'(s)}{R^{3-\alpha}}|\gamma'(\st+s)|\,\d \st\,\d z,
\]
for all $s\in\TT^1$.
The remaining term can be estimated as follows:
\begin{multline*}
\bigO{\varepsilon^{-2}}\left|\int\int_{\Sigma(\tau,\varepsilon)}\int_{-1/2}^{1/2}\frac{z\x\gamma'(s)}{R^{3-\alpha}}|\gamma'(\st+s)|\,\d \st\,\d z\,\d s\right|\\
\leq\bigO{\varepsilon^{-2}}\int\int_{\Sigma(\tau,\varepsilon)}|z|^{\alpha-1}\,\d z\,\d s\\
\leq \bigO{\varepsilon^{-2}} \int_{\tau-s\leq C_1\varepsilon/\inf\|\gamma'\|}\varepsilon^{\alpha+1}\,\d s\leq \bigO{\varepsilon^\alpha},
\end{multline*} 
where we have used that 
\[
\int \frac{1}{R^{3-\alpha}} \d \st=\int_{-1/2}^{-1/2} \frac{1}{\big||z|^2 + |\gamma'(\st)|^2 \st^2|\big|^{(3-\alpha)/2}} \d \st \leq C |z|^{\alpha-2},
\]
as can be seen by a simple scaling argument.

\subsection{Binormal term}
We next consider the binormal term $-\frac{\st^2}{2R^{3-\alpha}}\gamma'(s)\x\gamma''(s)$. Elementary calculations yield
\[
\frac{1-|z|^\alpha}{\alpha}\lesssim\int_{-1/2}^{1/2}\frac{\st^2|\gamma'(\st+s)|}{R^{3-\alpha}}\,\d \st\lesssim  \frac{1}{\alpha},
\]
hence 
\begin{multline*}
\lim_{\varepsilon\to 0}\mathcal{J}_\varepsilon \left(\int_{-1/2}^{1/2}\frac{-\st^2}{2R^{3-\alpha}}\gamma'(s)\x\gamma''(s)\,\d \st\right)(\tau)= -\frac{C(\alpha,\tau)}{\alpha}\gamma'(\tau)\x\gamma''(\tau)\\
=-\frac{C(\alpha,\tau)}{\alpha}|\gamma'(\tau)|^3\kappa(\tau) B(\tau),
\end{multline*}
where $C>0$ can be bounded above and below, independent of $\alpha$ and $\tau$.  Here $\kappa(\tau)$ and $B(\tau)$ denote the curvature and binormal to the curve $\gamma$ at $\gamma(\tau)$.
\subsection{Remaining terms}
More generally, for $k\geq1$, $m,\,n\geq 0$ we have
\begin{multline*}
\int_{-1/2}^{1/2}\frac{\st^kx^m}{R^{n-\alpha}}\,\d \st\leq \bigO{|x|^{k+m+\alpha-n+1}}\int_0^{1/|x|}\frac{\st^k}{(1+|\gamma'(s)|^2\st^2)^{(n-\alpha)/2}}\,\d \st\\
\leq  \bigO{|x|^{k+m+\alpha-n+1}}+  \frac{1}{\alpha+k+1-n}\bigO{|x|^{m}}.
\end{multline*}
Now for each but the leading term in \re{eq_expn_unmol} we have $k+m\geq n-1$. Thus each such term  gives a contribution to $u_\varepsilon(\gamma(\tau))$ of
\[\bigO{\varepsilon^\alpha} +  \frac{1}{\alpha+k+1-n}\bigO{\varepsilon^{m}}.\]

Noting that $k\geq n$ or $m\geq 1$ for all except for the binormal term, combining all of the estimates above yields
\begin{equation}\label{eqFinal}
\lim_{\varepsilon\to 0}u_\varepsilon(\gamma(\tau))=\frac{C(\alpha,\tau)}{\alpha}\kappa(\tau)B(\tau)+ w(\alpha,\tau),
\end{equation}
where $|w(\alpha,\tau)|$ is bounded independent of $\alpha,\tau$. 

Combining  the estimates from Section 3  the proof of Theorem \ref{thm1} is complete.

\section{Acknowledgements}
Both authors are partially supported by European Research Council, ERC Consolidator Grant no. 616797.


\end{document}